\documentclass[a4paper,10pt]{article}

\usepackage{amsmath,amssymb,amsfonts}
\usepackage{ifthen}
\usepackage{hyperref}
\usepackage[amsmath,thmmarks,hyperref]{ntheorem}

\theoremstyle{plain}
\theoremseparator{.}
\newtheorem{theorem}{Theorem}
\newtheorem{proposition}[theorem]{Proposition}
\newtheorem{lemma}[theorem]{Lemma}

\newtheorem{assumption}[theorem]{Assumption}

\theorembodyfont{}
\theoremsymbol{\ensuremath{{}_\blacksquare}}

\newtheorem{remark}[theorem]{Remark}

\theoremsymbol{\ensuremath{\square}}
\theoremheaderfont{\itshape}
\newtheorem*{proof}{Proof}

\newcommand{\req}[1]{\eqref{eq:#1}}
\DeclareMathOperator{\Span}{span}

\DeclareMathOperator{\sgn}{sgn}

\DeclareMathOperator*{\argmin}{arg\,min}
\DeclareMathOperator{\range}{range}
\DeclareMathOperator{\domain}{dom}

\newcommand{\field}[1]{\ensuremath{\mathbb{#1}}}
\newcommand{\R}{\field{R}}
\newcommand{\N}{\field{N}}

\newcommand{\inner}[3][n]{\SwitchBracketsizeLeft{#1}\LeftBracketSize\langle#2,#3\SwitchBracketsizeRight{#1}\RightBracketSize\rangle}
\newcommand{\abs}[2][n]{\SwitchBracketsizeLeft{#1}\LeftBracketSize\lvert#2\SwitchBracketsizeRight{#1}\RightBracketSize\rvert}
\newcommand{\norm}[2][n]{\SwitchBracketsizeLeft{#1}\LeftBracketSize\lVert#2\SwitchBracketsizeRight{#1}\RightBracketSize\rVert}
\newcommand{\set}[3][b]{\SwitchBracketsizeLeft{#1}\LeftBracketSize\{#2:#3\SwitchBracketsizeRight{#1}\RightBracketSize\}}

\newcommand{\NextScriptStyle}[1]{{\scriptstyle{#1}}}
\newcommand{\NextScriptScriptStyle}[1]{{\scriptscriptstyle{#1}}}
\newcommand{\NextTextStyle}[1]{{\textstyle{#1}}}
\newcommand{\NextDisplayStyle}[1]{{\displaystyle{#1}}}
\newcommand{\SwitchBracketsizeLeft}[1]{
  \ifthenelse{\equal{#1}{b}\OR\equal{#1}{big}}{\let\LeftBracketSize=\bigl}{
    \ifthenelse{\equal{#1}{B}\OR\equal{#1}{Big}}{\let\LeftBracketSize=\Bigl}{
      \ifthenelse{\equal{#1}{g}\OR\equal{#1}{bigg}}{\let\LeftBracketSize=\biggl}{
     \ifthenelse{\equal{#1}{G}\OR\equal{#1}{Bigg}}{\let\LeftBracketSize=\Biggl}{
      \ifthenelse{\equal{#1}{s}\OR\equal{#1}{small}}{\let\LeftBracketSize=\NextScriptStyle}{
        \ifthenelse{\equal{#1}{ss}}{\let\LeftBracketSize=\NextScriptScriptStyle}{
          \ifthenelse{\equal{#1}{t}\OR\equal{#1}{text}}{\let\LeftBracketSize=\NextTextStyle}{
         \ifthenelse{\equal{#1}{d}\OR\equal{#1}{display}}{\let\LeftBracketSize=\NextDisplayStyle}{
          \ifthenelse{\equal{#1}{a}\OR\equal{#1}{auto}}{\let\LeftBracketSize=\left}{
            \let\LeftBracketSize=\relax}}}}}}}}}}
\newcommand{\SwitchBracketsizeRight}[1]{
  \ifthenelse{\equal{#1}{b}\OR\equal{#1}{big}}{\let\RightBracketSize=\bigr}{
    \ifthenelse{\equal{#1}{B}\OR\equal{#1}{Big}}{\let\RightBracketSize=\Bigr}{
      \ifthenelse{\equal{#1}{g}\OR\equal{#1}{bigg}}{\let\RightBracketSize=\biggr}{
     \ifthenelse{\equal{#1}{G}\OR\equal{#1}{Bigg}}{\let\RightBracketSize=\Biggr}{
      \ifthenelse{\equal{#1}{s}\OR\equal{#1}{small}}{\let\RightBracketSize=\NextScriptStyle}{
        \ifthenelse{\equal{#1}{ss}}{\let\RightBracketSize=\NextScriptScriptStyle}{
          \ifthenelse{\equal{#1}{t}\OR\equal{#1}{text}}{\let\RightBracketSize=\NextTextStyle}{
         \ifthenelse{\equal{#1}{d}\OR\equal{#1}{display}}{\let\RightBracketSize=\NextDisplayStyle}{
          \ifthenelse{\equal{#1}{a}\OR\equal{#1}{auto}}{\let\RightBracketSize=\right}{
            \let\RightBracketSize=\relax}}}}}}}}}}

\makeatletter
\newcommand{\logmessage}[1]{\@latex@warning{#1}}

\newcommand{\ignore}{\logmessage{Text ignored}\@gobble}
\makeatother

\title{Sparse Regularization with $l^q$ Penalty Term}
\author{Markus Grasmair$^1$ \and Markus Haltmeier$^1$ \and Otmar Scherzer$^{1,2}$\\
\\
{\normalsize
\begin{tabular}{ccc}
\hbox to 0pt{\hss ${}^1$}Department of Mathematics && \hbox to 0pt{\hss ${}^2$}Radon Institute of Computational\\
  University of Innsbruck &&  and Applied Mathematics\\
  Technikerstr.~21a && Altenberger Str.~69\\
  6020 Innsbruck, Austria && 4040 Linz, Austria\\
\end{tabular}}}
\date{\today}

\begin{document}

\maketitle

\begin{abstract}
We consider the stable approximation of sparse solutions to non-linear operator equations by means
of Tikhonov regularization with a subqua\-dratic penalty term.
Imposing certain assumptions, which for a linear operator are equivalent
to the standard range condition, we derive the usual convergence rate $O(\sqrt{\delta})$
of the regularized solutions in dependence of the noise level $\delta$.
Particular emphasis lies on the case, where the true solution is known
to have a sparse representation in a given basis. In this case, if the differential of the operator
satisfies a certain injectivity condition, we can show that the actual convergence rate
improves up to $O(\delta)$.
\medskip

\end{abstract}

\section{Introduction}

A widely used technique for the
approximate solution of an ill-posed, possibly non-linear operator equation
\begin{equation}
\label{eq:operator}
F(u) = v
\end{equation}
on a Hilbert space $U$
is Tikhonov regularization, which can be formulated as minimization of
the functional
\[
{\cal T}(u) = \norm{F(u)-v}^2 + \alpha {\cal R}(u)\;.
\]
The first term ensures that the minimizer $u_\alpha$ will indeed approximately solve
the equation, while the second term stabilizes the process of inverting $F$
and forces $u_\alpha$ to satisfy certain regularity properties incorporated into ${\cal R}$.
Originally, Tikhonov applied this method to the stable solution of the Fredholm equation.
Requiring differentiability of $u_\alpha$,
he used the square of a higher order weighted Sobolev norm as penalty
term \cite{Tik63b_russ,Tik63a_russ}.

Recently, the focus has shifted from the postulation of differentiability properties
to sparsity constraints \cite{DauDefDeM04,BreLor08b,CanRomTao06,ComPes07,ComWaj05,DauForLor08,Don06,FigNowWri07,GriLor08,RamTes06,Tro06}.
Here, one requires the expansion of $u_\alpha$ with respect to
some given orthonormal basis $(\phi_i)_{i\in\N}$ of $U$
to be sparse in the sense that only finitely
many coefficients are different from zero.
This can be achieved with regularization functionals
\begin{equation}\label{eq:R}
    {\cal R}(u) = \sum_{i\in\N} \abs{\inner{\phi_i}{u}}^q\,,\qquad 0 \le q \le 2\;.
\end{equation}

In fact, sparsity of the solution is not necessarily guaranteed for $q > 1$.
The lack of convexity of ${\cal R}$, however, makes a choice $q < 1$
inconvenient both for theoretical analysis and the actual
computation of a minimizer. On the other hand, the assumption $q \leq 2$ is used to obtain
coercivity of the regularization functional, which in turn implies the existence of minimizers of ${\cal T}$.
For these reasons we only consider the case $1 \le q \le 2$.
\medskip

We concentrate our analysis on the well-posedness of the regularization method
and the derivation of convergence rates.
For that purpose we assume that only noisy data $v^\delta$
is given, which satisfies $\norm{v^\delta - v} \le \delta$.
We denote by $u_\alpha^\delta$ the minimizer of the regularization functional
with noisy data $v^\delta$ and regularization parameter $\alpha$,
and by $u^\dagger$ an ${\cal R}_q$-minimizing solution of $F(u) =  v$.
Then the question is how the distance $\norm{u_\alpha^\delta-u^\dagger}$
depends on the noise level $\delta$ and the regularization parameter $\alpha$.

Dismissing for the moment the assumption of sparsity, we derive
for a parameter choice $\alpha \sim \delta$ a convergence rate $\norm{u_\alpha^\delta - u^\dagger} = O(\sqrt{\delta})$
provided $1 < q \leq 2$ and a source condition is satisfied (see Proposition~\ref{prop:hkps}).
In the linear case this condition is the usual range condition
$\partial{\cal R}(u^\dagger) \cap \range(F^*) \neq \emptyset$,
where $F^*$ denotes the adjoint of the operator $F$ (see Proposition~\ref{prop:linear-aqu}).
Similar results have been derived recently \cite{Lor08,Ram08}.
In the non-linear case we impose a different assumption, which for sparsity regularization
generalizes common source conditions involving the Bregman distance \cite{HofKalPoeSch07,ResSch06,SchGraGroHalLen07}.

If, furthermore, the solution $u^\dagger$ of the operator equation is
known to be sparse, then the convergence rates of the regularized solutions to $u^\dagger$
can be shown to be $O(\delta^{1/q})$ where $1 \leq q \leq 2$ is the exponent in the regularization
term~\req{R} (see Theorems~\ref{th:rates_sparse_q} and~\ref{th:rates_sparse_1}).
To that end we require the derivative of $F$ at $u^\dagger$ to be invertible
on certain finite dimensional subspaces, a condition introduced in \cite{BreLor08}
for linear operators as `finite basis injectivity property'.
This improved convergence rate provides a theoretical justification for the
usage of subquadratic penalty terms for regularization with sparsity constraints.
\smallskip

Our results reveal a fundamental difference between quadratic and non-quadratic Tikhonov
regularization. Neubauer \cite{Neu97} has derived a saturation result
for quadratic regularization in a Hilbert space setting with a linear operator $F$.
He has shown that, apart from the trivial case $u^\dagger = 0$,
the convergence rates cannot be better than $O(\delta^{2/3})$.
The present article shows that this rate can be beaten by sparse regularization
when applied to the recovery of sparse data.

\section{Notational Preliminaries}

All along this paper we assume that $V$ is a reflexive Banach space
and $U$ is a Hilbert space in which a frame $(\phi_i)_{i\in\N} \subset U$ is given.
That is, there exist $0 < C_1 \leq C_2 < \infty$ such that
\[
C_1 \norm{u}^2 \le \sum_{i\in\N} \abs{\inner{\phi_i}{u}}^2 \le C_2\norm{u}^2
\qquad
\text{ for every } u \in U\;.
\]
The operator $F\colon \domain(F) \subseteq U \to V$ is assumed to be weakly sequentially
closed and $\domain(F) \cap \domain({\cal R}_q) \neq \emptyset$. Examples for weakly sequentially
closed operators are linear bounded operators restricted to convex domains, which naturally
arise for instance in image restoration problems or tomographic applications
\cite{EngHanNeu96,SchGraGroHalLen07}. Truly nonlinear operators arise in schlieren
imaging \cite{SchGraGroHalLen07} or simultaneous activity and attenuation reconstruction
in emission tomography \cite{Dic99}. See also \cite{BonBreLorMaa07,RamTes06} for the application
of sparsity constraints to inverse problems.
\smallskip

We define the regularization functional ${\cal R}_q \colon U \to \R \cup\{\infty\}$ by
\[
{\cal R}_q(u) := \sum_{i \in \N} w_i \abs{\inner{\phi_i}{u}}^q\,,
\]
where $1 \le q \le 2$ and there exists $w_{\min} > 0$ such that $w_i \ge w_{\min}$
for all $i \in \N$.
Note that ${\cal R}_q$ is convex and weakly lower semi-continuous as
the sum of non-negative convex and weakly continuous functionals.

The subdifferential of ${\cal R}_q$ at $u$ is denoted by $\partial {\cal R}_q (u) \subset U$.
If  $q>1$, then $\partial {\cal R}_q (u)$ is at most single valued and is identified with its
single element.

For the approximate solution of the operator equation $F(u) = v$
we consider the minimization of the regularization functional
\[
{\cal T}_{\alpha,v}^{p,q}(u) :=
\begin{cases}
  \norm[b]{ F(u)-v }^p + \alpha{\cal R}_q(u)\,,
  &   \text{if } u \in \domain(F) \cap \domain({\cal R}_q)\,,\\
  +   \infty\,,&\text{if } u \not\in\domain(F)\cap \domain({\cal R}_q)\,,
\end{cases}
\]
with some $\alpha > 0$ and $p\geq 1$.
\smallskip

In order to prove convergence rates results we impose an additional assumption concerning
the interaction of $F$ and ${\cal R}_q$ in a neighborhood of an ${\cal R}_q$-minimizing
solution of $F(u)=v$.
Here $u^\dagger \in U$ is called ${\cal R}_q$-minimizing solution, if
$F(u^\dagger) = v$ and
\[
    {\cal R}_q(u^\dagger) = \min\set{{\cal R}_q(u)}{F(u) = v}\;.
\]

\begin{assumption}\label{as:ratesr}
The equation $F(u) = v$ has an ${\cal R}_q$-minimizing solution $u^\dagger$
and there exist $\beta_1$, $\beta_2 > 0$, $r>0$, $\sigma > 0$, and $\rho > {\cal R}_q(u^\dagger)$ such that
\begin{equation}\label{eq:ratesr}
  {\cal R}_q(u)-{\cal R}_q(u^\dagger) \ge \beta_1 \norm{u-u^\dagger}^r - \beta_2 \norm{F(u)-F(u^\dagger)}
\end{equation}
for all $u \in \domain(F)$ satisfying ${\cal R}_q(u) < \rho$ and\/ $\norm{F(u)-F(u^\dagger)} < \sigma$.
\end{assumption}

In Section~\ref{se:source} below we show that Assumption~\ref{as:ratesr} with $r = 2$ follows
from the standard conditions stated in general convergence rates results in a Banach
space setting \cite{BurOsh04,HofKalPoeSch07,ResSch06}, which in turn generalize the standard conditions
in a Hilbert space setting \cite{EngHanNeu96,EngKunNeu89}.
Moreover, the assumption is equivalent to the
standard source condition $\partial{\cal R}(u^\dagger) \cap \range(F^*) \neq \emptyset$
in the particular case of a linear and bounded operator $F$ (see Proposition~\ref{prop:linear-aqu}).

\section{Well-Posedness and Convergence Rates}

In this section we prove the well-posedness of the regularization method.
By this we mean that minimizers $u_\alpha^\delta$ of the regularization functional
${\cal T}_{\alpha,v^\delta}^{p,q}$ exist for every $\alpha > 0$, continuously
depend on the data $v^\delta$, and converge to a solution of $F(u) = v$ as the
noise level approaches zero, provided the regularization parameter $\alpha$ is chosen
appropriately.

These results are analogous to results obtained for standard quadratic Ti\-kho\-nov
regularization in Hilbert spaces (see e.g.~\cite{EngHanNeu96}). Also the mathematical techniques
employed in the proofs of existence, weak stability, and convergence are similar.
Some extra work is needed, however, for the passage from weak stability and convergence
to stability and convergence with respect to ${\cal R}_q$.

\begin{lemma}\label{le:weakstrong}
Let $1\leq q\leq 2$. Assume that $(u_k)_{k\in\N}\subset U$\! weakly converges to $u\in U$\! and that
${\cal R}_q(u_k)$ converges to ${\cal R}_q(u)$. Then ${\cal R}_q(u_k-u) \to 0$.
\end{lemma}

\begin{proof}
The assumption ${\cal R}_q (u_k) \to {\cal R}(u)$ implies that
\begin{multline*}
\limsup_k {\cal R}_q(u_k-u) \\
\begin{aligned}
&= \limsup_k \Bigl[2\bigl({\cal R}_q(u_k) + {\cal R}_q(u)\bigr) -
  2\bigl({\cal R}_q(u_k)+{\cal R}_q(u)\bigr)+{\cal R}_q(u_k-u)\Bigr]\\
&= 4\,{\cal R}_q(u)
  - \liminf_k \sum_{i\in\N} w_i \Bigl[2\,\abs{\inner{\phi_i}{u_k}}^q+2\,\abs{\inner{\phi_i}{u}}^q
  - \abs{\inner{\phi_i}{u_k-u}}^q\Bigr]\;.
\end{aligned}
\end{multline*}
Using Fatou's Lemma we obtain that
\begin{multline*}
- \liminf_k \sum_{i\in\N} w_i \Bigl[2\,\abs{\inner{\phi_i}{u_k}}^q+2\,\abs{\inner{\phi_i}{u}}^q
  - \abs{\inner{\phi_i}{u_k-u}}^q\Bigr]\\
\le - \sum_{i\in\N}\liminf_k  w_i \Bigl[2\,\abs{\inner{\phi_i}{u_k}}^q+2\,\abs{\inner{\phi_i}{u}}^q
  - \abs{\inner{\phi_i}{u_k-u}}^q\Bigr]\;.
\end{multline*}
Now, the weak convergence of $(u_k)_{k\in\N}$ shows that $\inner{\phi_i}{u_k} \to \inner{\phi_i}{u}$
for all $i \in \N$.
Therefore it follows that
\[
- \sum_{i\in\N} \liminf_k w_i \Bigl[2\,\abs{\inner{\phi_i}{u_k}}^q+2\,\abs{\inner{\phi_i}{u}}^q
  - \abs{\inner{\phi_i}{u_k-u}}^q\Bigr]
= -4\sum_{i\in\N} w_i\abs{\inner{\phi_i}{u}}^q\;.
\]
Combining the above inequality and equalities we see that
\[
\limsup_k {\cal R}_q(u_k-u) \leq 4\,{\cal R}_q(u) - 4\sum_{i\in\N} w_i\abs{\inner{\phi_i}{u}}^q = 0
\]
or, equivalently, that ${\cal R}_q(u_k-u) \to 0$.
\end{proof}

\begin{remark}
Convergence with respect to ${\cal R}_q$ implies convergence with respect to the norm, which is an
easy consequence of the inequality
\begin{equation}\label{eq:ineq}
\Bigl( \sum_{i\in\N} \abs{c_i}^t \Bigr)^{1/t}
\le
\Bigl( \sum_{i\in\N} \abs{c_i}^s \Bigr)^{1/s} =: \abs{c}_s
\end{equation}
for $c = (c_i)_{i\in\N} \in \R^\N$ and $0 < s \le t < \infty$.
The inequality~\req{ineq} easily follows for $0 < \abs{c}_s < \infty$
from the inequality
\[
\sum_{i\in\N}\biggl( \frac{\abs{c_i}}{\abs{c}_s}\biggr)^t
\le \sum_{i\in\N} \biggl(\frac{\abs{c_i}}{\abs{c}_s}\biggr)^s  = 1 \;.
\]
In particular, this shows that
\begin{equation}\label{eq:ineq2}
{\cal R}_q(u)
\ge w_{\min} \sum_{i\in\N} \abs{\inner{\phi_i}{u}}^q
\ge w_{\min} \Bigl(\sum_{i\in\N} \abs{\inner{\phi_i}{u}}^2\Bigr)^{q/2}
\ge w_{\min} C_1^{q/2} \norm{u}^q
\end{equation}
for every $u \in U$.
Therefore, Lemma~\ref{le:weakstrong} implies \cite[Lemma 4.3]{DauDefDeM04},
where the authors show convergence of the sequence $(u_k)_{k\in\N}$ with respect to the norm.

Another immediate consequence of~\req{ineq2} is the weak coercivity of
the functional ${\cal R}_q$.
\end{remark}

\begin{lemma}\label{le:weakconvergence}
Let $(u_k)_{k \in \N} \subset \domain(F)$ and $(v_k)_{k\in\N}
\subset V$. Assume that the sequence $(v_k)_{k\in\N}$ is bounded in
$V$\! and that there exist $\alpha > 0$ and $M > 0$ such that ${\cal
T}_{\alpha,v_k}^{p,q}(u_k) < M$ for all $k \in \N$. Then there exist
$u \in \domain(F)$ and a subsequence $(u_{k_j})_{j\in\N}$ such that
$u_{k_j} \rightharpoonup u$ and $F(u_{k_j}) \rightharpoonup F(u)$.
\end{lemma}

\begin{proof}
The coercivity of ${\cal R}_q$ and the estimate $ {\cal T}_{\alpha,v_k}^{p,q}(u_k) \ge \alpha {\cal R}_q(u_k)$
imply that the sequence $(u_k)_{k\in\N}$ is bounded in $U$.
Similarly, since $(v_k)_{k\in\N}$ is bounded,
also the sequence $\bigl(F(u_k)\bigr)_{k\in\N}$ is bounded in $V$.
Therefore there exist a subsequence $(u_{k_j})_{j\in\N}$ and $u \in U$, $y \in V$,
such that $(u_{k_j})_{j\in\N}$ weakly converges to $u$ and $\bigl(F(u_{k_j})\bigr)_{j\in\N}$
weakly converges to $y$.
Since $F$ is weakly sequentially closed, it follows that $u \in \domain(F)$
and $F(u) = y$.
\end{proof}

The ideas of the following proofs are based on~\cite[Section 3]{HofKalPoeSch07}.
Still, we provide short proofs, since our assumptions are
slightly different from~\cite{HofKalPoeSch07}, where weak continuity
of the operator $F$ is assumed.

\begin{proposition}[Existence]
For every $v^\delta \in V$ the functional ${\cal T}_{\alpha,v^\delta}^{p,q}$ has a minimizer in $U$.
\end{proposition}

\begin{proof}
Let $(u_k)_{k\in\N}$ satisfy
\[
\lim_{k\to \infty} {\cal T}_{\alpha,v^\delta}^{p,q}(u_k)
= \inf \{ {\cal T}_{\alpha,v^\delta}^{p,q}(u): u \in U \}\;.
\]
Lemma~\ref{le:weakconvergence} shows that there exists a subsequence $(u_{k_j})_{j\in\N}$
weakly converging to some $u \in U$ such that $F(u_{k_j}) \rightharpoonup F(u)$.
Therefore the weak sequential lower semi-continuity of ${\cal T}_{\alpha,v^\delta}^{p,q}$
implies that $u$ is a minimizer of ${\cal T}_{\alpha,v^\delta}^{p,q}$.
\end{proof}

\begin{proposition}[Stability]
Let $(v_k)_{k\in\N}$ converge to $v^\delta \in V$\!
and let
\[
u_k \in \argmin\set{{\cal T}_{\alpha,v_k}^{p,q}(u)}{u\in U}\;.
\]
Then there exists a subsequence $(u_{k_j})_{j\in\N}$ and a
minimizer \ $u_\alpha^\delta$ of \ ${\cal T}_{\alpha,v^\delta}^{p,q}$ such that
${\cal R}_q(u_\alpha^\delta - u_{k_j}) \to 0$.
If the minimizer $u_\alpha^\delta$ is unique, then $(u_k)_{k\in\N}$ converges to
$u_\alpha^\delta$ with respect to ${\cal R}_q$.
\end{proposition}

\begin{proof}
From Lemma~\ref{le:weakconvergence} we obtain the existence of a subsequence $(u_{k_j})_{j\in\N}$
weakly converging to some $u \in \domain(F)$ such that $F(u_{k_j})\rightharpoonup F(u)$.
Since $v_k \to v^\delta$, it follows that
${\cal T}_{\alpha,v^\delta}^{p,q}(u) \le \liminf_j {\cal T}_{\alpha,v_{k_j}}^{p,q}(u_{k_j})$.

On the other hand, if $\tilde{u} \in \domain(F)$, then
\[
{\cal T}_{\alpha,v^\delta}^{p,q}(\tilde{u}) = \lim_k {\cal T}_{\alpha,v_k}^{p,q}(\tilde{u})
\ge \liminf_k {\cal T}_{\alpha,v_k}^{p,q}(u_k)\;.
\]
Thus $u = u_\alpha^\delta$ is a minimizer of ${\cal T}_{\alpha,v^\delta}^{p,q}$.

Now note that also ${\cal T}_{\alpha,v^\delta}^{p,q}(u_{k_j}) \to {\cal T}_{\alpha,v^\delta}^{p,q}(u)$.
Since both $\norm{\cdot}^p$ and ${\cal R}_q$ are weakly sequentially lower semi-continuous,
this implies that ${\cal R}_q(u_{k_j}) \to {\cal R}_q(u)$.
Using Lemma~\ref{le:weakstrong}, we therefore obtain the convergence of the sequence
$(u_{k_j})_{j\in\N}$ with respect to ${\cal R}_q$.

In case the minimizer $u_\alpha^\delta$ is unique, the convergence of the original
sequence $(u_k)_{k\in\N}$ to $u_\alpha^\delta$ follows from a subsequence argument.
\end{proof}

\begin{proposition}[Convergence]
Assume that the operator equation $F(u) = v$ attains a solution in $\domain({\cal R}_q)$
and that $\alpha\colon \R_{> 0} \to \R_{> 0}$ satisfies
\[
\alpha(\delta) \to 0
\qquad
\text{and}
\qquad
\frac{\delta^p}{\alpha(\delta)} \to 0
\qquad
\text{as }
\delta \to 0\;.
\]
Let $\delta_k \to 0$ and let $v_k \in V$ satisfy $\norm{v_k-v} \le \delta_k$.
Moreover, let $\alpha_k = \alpha(\delta_k)$ and
\[
u_k \in \argmin\set{{\cal T}_{\alpha_k,v_k}^{p,q}\!(u)}{u\in U}\;.
\]

Then there exist an ${\cal R}_q$-minimizing solution $u^\dagger$ of $F(u)=v$
and a subsequence $(u_{k_j})_{j\in\N}$ with
${\cal R}_q(u^\dagger - u_{k_j}) \to 0$.
If the ${\cal R}_q$-minimizing solution is unique,
then $(u_k)_{k\in\N}$ converges to $u^\dagger$ with respect to ${\cal R}_q$.
\end{proposition}

\begin{proof}
Let $\tilde{u} \in \domain({\cal R}_q)$ be any solution of $F(u) = v$.
The definition of $u_k$ implies that
\[
\norm{F(u_k)-v_k}^p + \alpha_k{\cal R}_q(u_k)
\le \norm{F(\tilde{u})-v_k}^p + \alpha_k {\cal R}_q(\tilde{u})
\le \delta_k^p + \alpha_k{\cal R}_q(\tilde{u})\;.
\]
In particular $\norm{F(u_k)-v_k} \to 0$ and
\begin{equation}\label{eq:convergence:limsup}
\limsup_k {\cal R}_q(u_k)
\le {\cal R}_q(\tilde{u}) + \limsup_k \frac{\delta^p_k}{\alpha_k}
= {\cal R}_q(\tilde{u})\;.
\end{equation}

This shows that there exists $M > 0$ such that ${\cal T}_{\alpha_1,v_k}^{p,q}(u_k) \le M$
for all $k \in \N$.
Thus Lemma~\ref{le:weakconvergence} yields a subsequence
$(u_{k_j})_{j\in\N}$ weakly converging to some $u^\dagger \in \domain(F)$
such that $F(u_{k_j}) \rightharpoonup F(u^\dagger)$.
Since $\norm{F(u_{k_j})-v} \le \norm{F(u_{k_j})-v_{k_j}} + \norm{v_{k_j}-v} \to 0$, it follows that
$F(u^\dagger) = v$.

The weak sequential lower semi-continuity of ${\cal R}_q$ implies that
${\cal R}_q(u^\dagger) \le \liminf_j {\cal R}_q(u_{k_j})$.
Since~\req{convergence:limsup} holds for every $\tilde{u} \in \domain({\cal R}_q)$
satisfying $F(\tilde{u}) = v$, it follows that $u^\dagger$ is an ${\cal R}_q$-minimizing
solution of $F(u) = v$ and that ${\cal R}_q(u_{k_j}) \to {\cal R}_q(u^\dagger)$.
Lemma~\ref{le:weakstrong} now shows that $(u_{k_j})_{j\in\N}$
converges to $u^\dagger$ with respect to ${\cal R}_q$.

Again, the convergence of the original sequence $(u_k)_{k\in\N}$ to $u^\dagger$
follows from a subsequence argument, if the ${\cal R}_q$-minimizing solution $u^\dagger$ is unique.
\end{proof}

\smallskip
In the following we write $\alpha \sim \delta^s$ for
$\alpha\colon (0,\infty) \to (0,\infty)$ and $s > 0$, if there exist constants
$C \ge c > 0$ and $\delta_0>0$,
such that $c \delta^s \leq \alpha(\delta) \leq C \delta^s$ for every
$0<\delta<\delta_0$.

\smallskip
For the next result on convergence rates recall the definition of the exponent $r$
in Assumption~\ref{as:ratesr}.

\begin{proposition}[Convergence Rates]
\label{pr:con_rates}
Let Assumption~\ref{as:ratesr} hold.
Assume that $v^\delta \in V$\! satisfies $\norm{v^\delta-v} \le \delta$
and $u_\alpha^\delta \in \argmin\set{{\cal T}_{\alpha,v^\delta}^{p,q}(u)}{u\in U}$.
For $\alpha$ and $\delta$ sufficiently small we obtain the following estimates:

If $p=1$ and $\alpha\beta_2 < 1$, then
\[
\norm{u_\alpha^\delta-u^\dagger}^r \le \frac{(1+\alpha\beta_2)\,\delta}{\alpha\beta_1}\,,
\qquad
\norm{F(u_\alpha^\delta)-v^\delta} \le \frac{(1+\alpha\beta_2)\,\delta}{1-\alpha\beta_2}\;.
\]

If $p>1$, then
\[
\begin{aligned}
\norm{u_\alpha^\delta-u^\dagger}^r &\le \frac{\delta^p + \alpha\beta_2\delta + (\alpha\beta_2)^{p_*}/p_*}{\alpha\beta_1}\,,\\
\norm{F(u_\alpha^\delta)-v^\delta}^p &\le p_* \delta^p + p_*\alpha\beta_2\delta + (\alpha\beta_2)^{p_*}\;.
\end{aligned}
\]
Here, $p_*$ is the conjugate of $p$ defined by $1/p_*+1/p = 1$.

In particular, if $\alpha \sim \delta^{p-1}$, then $\norm{u_\alpha^\delta - u^\dagger} = O(\delta^{1/r})$.
\end{proposition}

\begin{proof}
Since $u_\alpha^\delta$ minimizes ${\cal T}_{\alpha,v^\delta}^{p,q}$, the inequality
\[
\norm{F(u_\alpha^\delta)-v^\delta}^p + \alpha{\cal R}_q(u_\alpha^\delta)
\le \norm{F(u^\dagger)-v^\delta}^p + \alpha{\cal R}_q(u^\dagger)
\]
holds. Assumption~\ref{as:ratesr} and the fact that $F(u^\dagger) = v$ therefore
imply that
\[
\begin{aligned}
\delta^p &
\ge \norm{F(u_\alpha^\delta)-v^\delta}^p + \alpha\bigl({\cal R}_q(u_\alpha^\delta) - {\cal R}_q(u^\dagger)\bigr)\\
&\ge \norm{F(u_\alpha^\delta)-v^\delta}^p + \alpha\,\beta_1\,\norm{u_\alpha^\delta-u^\dagger}^r
    - \alpha\,\beta_2\,\norm{F(u_\alpha^\delta)-F(u^\dagger)}\\
&\ge \norm{F(u_\alpha^\delta)-v^\delta}^p + \alpha\,\beta_1\,\norm{u_\alpha^\delta-u^\dagger}^r
   - \alpha\,\beta_2\,\norm{F(u_\alpha^\delta) - v^\delta} - \alpha\beta_2\delta\;.
\end{aligned}
\]
This shows the assertion in the case $p=1$.

If $p > 1$, we apply Young's inequality $ab \le a^p/p + b^{p_*}/p_*$ with
$a = \norm{F(u_\alpha^\delta)-v^\delta}$ and $b = \alpha\beta_2$.
Then again the assertion follows.
\end{proof}

\begin{remark}
 Proposition \ref{pr:con_rates} shows that sparsity regularization is an \emph{exact meth\-od} for $p=1$,
that is, it yields exact solutions $u^\dagger$ for noise free data and $\alpha < 1/\beta_2$.
\end{remark}

\section{Relations to Source Conditions}\label{se:source}

We now investigate Assumption~\ref{as:ratesr} more closely and show
that it is indeed a generalization of commonly imposed source conditions
involving the Bregman distance defined by the functional ${\cal R}_q$
(see e.g.~\cite{BurOsh04,HofKalPoeSch07}).
The basis of these results is the following lemma, which relates the Bregman distance
to the squared norm on $U$ in case $q > 1$. This result is a consequence of a
special case of \cite[Lemma 2.7]{BonKazMaaSchoSchu08} (see also
\cite[Corollary 3.7]{ButIusZal03}).

From now on we assume that $(\phi_i)_{i\in \N}$ is an orthonormal basis.

\begin{lemma}\label{le:bregman_norm}
Let\/ $1 < q \le 2$. There exists a constant $c_q>0$ only depending on $q$ such that
\[
{\cal D}_B(\tilde{u},u) := {\cal R}_q(\tilde{u})-{\cal R}_q(u) - \inner{\partial{\cal R}_q(u)}{\tilde{u}-u}
\ge \frac{c_q\,\norm{\tilde{u}-u}^2}{3w_{\min} + 2{\cal R}_q(u)+{\cal R}_q(\tilde{u})}\,
\]
for all $\tilde{u}$, $u \in \domain({\cal R}_q)$ for which $\partial{\cal R}_q(u) \neq\emptyset$,
which is equivalent to the assumption that $\sum_{i\in\N} w_i^2 \, \abs{\inner{\phi_i}{u}}^{2(q-1)} < \infty$.
\end{lemma}

\begin{proof}
There exists $d_q > 0$ such that
\begin{equation}\label{eq:dq}
d_q \abs{a-b}^2
\le (\abs{a}^{2-q} + \abs{a-b}^{2-q}) \Bigl[\abs{b}^q - \abs{a}^q - q\abs{a}^{q-1}\sgn(a)\,(b-a)\Bigr]
\end{equation}
for all $a$, $b \in \R$ \cite[\S 5, Eq.~1]{Die75}.

Let $\tilde{u} \neq u \in \domain({\cal R}_q)$. Then
\[
\partial{\cal R}_q(u)
= \sum_{i\in\N} q\,w_i\,\abs{\inner{\phi_i}{u}}^{q-1}\sgn(\inner{\phi_i}{u})\,\phi_i
\]
provided that $\partial{\cal R}_q(u) \neq\emptyset$.
Applying~\req{dq}, we see that
\begin{multline}
\label{eq:bregman}
{\cal R}_q(\tilde{u}) - {\cal R}_q(u) - \inner{\partial {\cal R}_q}{\tilde{u}-u} \\
\begin{aligned}
&= \sum_{i\in\N} w_i\Bigl[ \abs{\inner{\phi_i}{\tilde{u}}}^q - \abs{\inner{\phi_i}{u}}^q
    - q\,\abs{\inner{\phi_i}{u}}^{q-1}\sgn(\inner{\phi_i}{u})\,\inner{\phi_i}{\tilde{u}-u}\Bigr]\\
&\ge d_q \sum_{i\in\N}
  \frac{w_i\,\abs{\inner{\phi_i}{\tilde{u}-u}}^2}
       {\abs{\inner{\phi_i}{u}}^{2-q}+\abs{\inner{\phi_i}{\tilde{u}-u}}^{2-q}}\\
&\ge \frac{d_q\,w_{\min}}{\max\set[n]{\abs{\inner{\phi_i}{u}}^{2-q}+\abs{\inner{\phi_i}{\tilde{u}-u}}^{2-q}}{i\in\N}}
       \sum_{i\in\N} \abs{\inner{\phi_i}{\tilde{u}-u}}^2\\
&\ge \frac{d_q\,w_{\min}}{\max\set[n]{2\abs{\inner{\phi_i}{u}}^{2-q}+\abs{\inner{\phi_i}{\tilde{u}}}^{2-q}}{i\in\N}}
       \sum_{i\in\N} \abs{\inner{\phi_i}{\tilde{u}-u}}^2\\
&\ge \frac{d_q\,w_{\min}}{3+\max\set[n]{2\abs{\inner{\phi_i}{u}}^q+\abs{\inner{\phi_i}{\tilde{u}}}^q}{i\in\N}}
       \sum_{i\in\N} \abs{\inner{\phi_i}{\tilde{u}-u}}^2\\
&\ge \frac{d_q\,w_{\min}^2}{3w_{\min}+2{\cal R}_q(u) + {\cal R}_q(\tilde{u})}
       \sum_{i\in\N} \abs{\inner{\phi_i}{\tilde{u}-u}}^2
       \;.
\end{aligned}
\end{multline}
Here, the third and second to last estimates follow from the inequalities
$(a+b)^{2-q} \leq a^{2-q} + b^{2-q}$ and $a^{q-2} \leq 1+a^q$ for $a$, $b\geq 0$.
Thus the assertion follows by setting $c_q := d_q\,w_{\min}^2$.
\end{proof}

\begin{proposition}\label{prop:linear-aqu}
Let $F$ be a bounded linear operator on $U$, $1 < q \le 2$, and $u^\dagger$ an ${\cal R}_q$-minimizing solution of $F(u) = v$.
Then Assumption~\ref{as:ratesr} with $r = 2$ is equivalent to the source condition
\begin{equation}\label{eq:source_linear}
\partial{\cal R}_q(u^\dagger) \in \range(F^*)\;.
\end{equation}
In particular, if $\alpha \sim \delta^{p-1}$, then $\norm{u_\alpha^\delta-u^\dagger} = O(\sqrt{\delta})$.
\end{proposition}

\begin{proof}
First assume that~\req{source_linear} holds.
The condition $\partial{\cal R}_q(u^\dagger) \in \range(F^*)$ implies the existence
of a constant $\hat{C} > 0$ such that
\begin{equation}\label{eq:contF}
    \abs{\inner{\partial{\cal R}_q(u^\dagger)}{u-u^\dagger}} \le \hat{C} \norm{F(u-u^\dagger)}
\end{equation}
for all $u \in U$.
Together with Lemma~\ref{le:bregman_norm} this yields the inequality
\begin{multline*}
{\cal R}_q(u)-{\cal R}_q(u^\dagger)
\ge \frac{c_q}{3w_{\min} + 2{\cal R}_q(u^\dagger)+{\cal R}_q(u)}\,\norm{u-u^\dagger}^2
    + \inner{\partial{\cal R}_q(u^\dagger)}{u-u^\dagger}\\
\ge \frac{c_q}{3w_{\min} + 2{\cal R}_q(u^\dagger)+{\cal R}_q(u)}\,\norm{u-u^\dagger}^2 - \hat{C}\norm{F(u-u^\dagger)}\;.
\end{multline*}
Thus, Assumption~\ref{as:ratesr} is satisfied if we choose
$r=2$, $\rho = {\cal R}_q(u^\dagger) +w_{\min}$, $\beta_1 = c_q/(4w_{\min}+3{\cal R}_q(u^\dagger))$,
and $\beta_2 = \hat{C}$.
\smallskip

In order to show the converse implication,
let Assumption~\ref{as:ratesr} be satisfied for $r=2$, that is, there exist $\beta_1$, $\beta_2 > 0$ such that
\[
\beta_1\norm{u-u^\dagger}^2 \le {\cal R}_q(u)-{\cal R}_q(u^\dagger) + \beta_2\norm{F(u-u^\dagger)}
\]
in a neighborhood of $u^\dagger$.
Both sides of this inequality are convex functions in the variable $u$ that agree for $u = u^\dagger$.
This implies that the subgradient at $u^\dagger$ of the left hand side, which equals zero, is contained in the
subgradient at $u^\dagger$ of the right hand side. In other words,
\[
0 \in \partial{\cal R}_q(u^\dagger) + \beta_2 F^* \partial\bigl(\norm{F(u-u^\dagger)}\bigr)\;.
\]
Consequently the source condition~\req{source_linear} holds.
\end{proof}

The following result states that the condition proposed in~\cite{HofKalPoeSch07}
for obtaining convergence rates in
the non-linear, non-smooth case also follows from Assumption~\ref{as:ratesr} with exponent $r=2$.

\begin{proposition}\label{prop:hkps}
Let\/ $1 < q <2$ and $u^\dagger$ an ${\cal R}_q$-minimizing solution of $F(u) = v$.
Assume that there exist $0 \le \gamma_1 < 1$, $\gamma_2 > 0$, and $\rho > {\cal R}_q(u^\dagger)$ such that
\begin{equation}\label{eq:ratesr_bregman}
\inner[b]{\partial {\cal R}_q(u^\dagger)}{u^\dagger - u} \le \gamma_1\,{\cal D}_B(u,u^\dagger) + \gamma_2\,\norm{F(u)-F(u^\dagger)}
\end{equation}
for all $u \in \domain(F)$ with ${\cal R}_q(u) < \rho$.
Then Assumption~\ref{as:ratesr} holds with $r=2$.
In particular, if $\alpha \sim \delta^{p-1}$, then $\norm{u_\alpha^\delta-u^\dagger} = O(\sqrt{\delta})$.
\end{proposition}

\begin{proof}
Using~\req{ratesr_bregman} and Lemma~\ref{le:bregman_norm} we obtain that
\begin{multline*}
\gamma_1 \bigl({\cal R}_q(u)-{\cal R}_q(u^\dagger)\bigr)
\ge -(1-\gamma_1) \inner[b]{\partial{\cal R}_q(u^\dagger)}{u-u^\dagger}-\gamma_2\,\norm{F(u)-F(u^\dagger)}\\
\ge \tilde{\beta}\norm{u-u^\dagger}^2 - \gamma_2\,\norm{F(u)-F(u^\dagger)}
    - (1-\gamma_1)\bigl({\cal R}_q(u)-{\cal R}_q(u^\dagger)\bigr)\,,
\end{multline*}
where $\tilde{\beta} := (1-\gamma_1)\,c_q/(3w_{\min}+2{\cal R}_q(u^\dagger) + \rho)$.
Thus Assumption~\ref{as:ratesr} follows with $\beta_1 = \tilde{\beta}/(1+2\gamma_1)$
and $\beta_2 = \gamma_2/(1+2\gamma_1)$.
\end{proof}

\section{Convergence Rates for Sparse Solutions}

We have seen above that appropriate source conditions imply convergence rates of
type $\sqrt{\delta}$.  These rates in fact can be improved considerably,
if the ${\cal R}_q$-minimizing solution $u^\dagger$ is sparse with respect to
$(\phi_i)_{i\in\N}$ in the sense that the set
\begin{equation*}
  J := \set{i \in \N}{\inner{u^\dagger}{\phi_i} \neq 0}
\end{equation*}
is finite.

\begin{assumption}\label{as:sparse_rates}
Assume that the following hold:
\begin{enumerate}
\item The operator equation $F(u) = v$ has an ${\cal R}_q$-minimizing
  solution $u^\dagger$ that is sparse with respect to $(\phi_i)_{i \in \N}$.
\item The operator $F$ is G\^ateaux differentiable at $u^\dagger$, and for
  every finite set $J \subset \N$ the restriction
  of its derivative $F'(u^\dagger)$ to $\Span\set{\phi_j}{j\in J}$ is injective.
\item There exist $\gamma_1$, $\gamma_2 > 0$, $\sigma > 0$, and $\rho > {\cal R}_q(u^\dagger)$ such that
\begin{equation}\label{eq:rates_sparse_q}
{\cal R}_q(u)-{\cal R}_q(u^\dagger)
\ge \gamma_1 \norm[b]{F(u)-F(u^\dagger)-F'(u^\dagger)(u-u^\dagger)} - \gamma_2 \norm[b]{F(u)-F(u^\dagger)}
\end{equation}
for all $u \in \domain(F)$ satisfying ${\cal R}_q(u) < \rho$ and\/ $\norm{F(u)-F(u^\dagger)} < \sigma$.
\end{enumerate}
\end{assumption}

We first derive a convergence rates result of
order $\delta^{1/q}$ for $q > 1$.

\begin{theorem}[$q>1$]\label{th:rates_sparse_q}
Let $1 < q \le 2$ and assume that Assumption~\ref{as:sparse_rates} holds.
Then for a parameter choice strategy $\alpha \sim \delta^{p-1}$
we obtain the convergence rate
\begin{equation*}
\norm{u_\alpha^\delta-u^\dagger} = O(\delta^{1/q})\;.
\end{equation*}
\end{theorem}

\begin{proof}
We verify Assumption~\ref{as:ratesr} with $r = q$ and appropriate constants
$\beta_1$, $\beta_2 > 0$.
Then the assertion follows from Proposition~\ref{pr:con_rates}.
\smallskip

Let therefore $u \in U$ satisfy ${\cal R}_q(u) < \rho$ and
$\norm{F(u)-F(u^\dagger)} < \sigma$.

Define $J := \set[n]{i\in\N}{\inner{u^\dagger}{\phi_i} \neq 0}$
and $W := \Span\set[n]{\phi_j}{j\in J}$.
Since $u^\dagger$ is sparse, the set $J$ is finite.
Therefore, the restriction of $F'(u^\dagger)$ to $W$ is injective,
which implies the existence of a constant $C > 0$ such that
\[
C\norm{F'(u^\dagger)\,w} \ge \norm{w}
\qquad
\text{ for all } w \in W\;.
\]

Now denote by $\pi_W$, $\pi_W^\perp\colon U \to U$
the projections
\[
\pi_W u := \sum_{j \in J} \inner{\phi_j}{u}\,\phi_j\,,
\qquad
\pi_W^\perp u := \sum_{j \not\in J} \inner{\phi_j}{u}\,\phi_j\;.
\]
Note that by assumption $\inner{\phi_j}{u^\dagger}=0$ for every $j \not \in J$,
which implies that $u^\dagger = \pi_W u^\dagger$ and $\pi_W^\perp u^\dagger = 0$.
By means of the inequality
\[
(a+b)^q \le 2^{q-1}(a^q + b^q) \le 2(a^q+b^q)
\qquad
\text{ for every } a,\, b > 0
\]
it therefore follows that
\begin{equation}\label{eq:rates_sparse_q:1}
\begin{aligned}
\norm{u-u^\dagger}^q
&\le 2\,\norm[b]{\pi_W(u-u^\dagger)}^q+2\norm[b]{\pi_W^\perp u}^q\\
&\le 2\, C^q\norm[b]{F'(u^\dagger)\bigl(\pi_W(u-u^\dagger)\bigr)}^q + 2\,\norm[b]{\pi_W^\perp u}^q\\
&\le 4\, C^q \norm[b]{F'(u^\dagger)(u-u^\dagger)}^q + 2\bigl(1+2 C^q\norm{F'(u^\dagger)}^q\bigr) \norm[b]{\pi_W^\perp u}^q\;.
\end{aligned}
\end{equation}

We now derive an estimate for $\norm{\pi_W^\perp u}^q$.
Using \req{ineq} we see that
\begin{equation}\label{eq:pibot1}
\norm[b]{\pi_W^\perp u}^q
=  \Bigl(\sum_{i\not\in J}\abs{\inner{\phi_i}{u}}^2\Bigr)^{q/2}
\le \sum_{i \not \in J} \abs{\inner{\phi_i}{u}}^q
\le w_{\min}^{-1}\sum_{i \not \in J} w_i \abs{\inner{\phi_i}{u}}^q\;.
\end{equation}
Since $q > 1$, the inequality
\[
\abs{\inner{\phi_i}{u}}^q - \abs{\inner{\phi_i}{u^\dagger}}^q
  - q \abs{\inner{\phi_i}{u^\dagger}}^{q-1}\sgn\bigl(\abs{\inner{\phi_i}{u^\dagger}}\bigr)\inner{\phi_i}{u-u^\dagger}
\geq 0
\]
holds for all $i \in \N$.
Consequently,
\begin{multline}\label{eq:pibot2}
\sum_{i\not\in J} w_i \abs{\inner{\phi_i}{u}}^q \\
\begin{aligned}
&= \sum_{i\not\in J} w_i \Bigl[\abs{\inner{\phi_i}{u}}^q - \abs{\inner{\phi_i}{u^\dagger}}^q
  - q \abs{\inner{\phi_i}{u^\dagger}}^{q-1}\sgn\bigl(\abs{\inner{\phi_i}{u^\dagger}}\bigr)
  \inner{\phi_i}{u-u^\dagger} \Bigr]\\
&\le \sum_{i\in \N} w_i \Bigl[ \abs{\inner{\phi_i}{u}}^q - \abs{\inner{\phi_i}{u^\dagger}}^q
  -  q \abs{\inner{\phi_i}{u^\dagger}}^{q-1}\sgn\bigl(\abs{\inner{\phi_i}{u^\dagger}}\bigr)\inner{\phi_i}{u-u^\dagger} \Bigr]\\
&= {\cal R}_q(u)-{\cal R}_q(u^\dagger) - \inner[b]{\partial{\cal R}_q(u^\dagger)}{u-u^\dagger}\;.
\end{aligned}
\end{multline}

From~\req{rates_sparse_q} we obtain by considering $u = u^\dagger+t\tilde u$, dividing by $t$, and passing
to the limit $t \to 0$ that
\begin{equation}\label{eq:subgradient}
\inner[b]{\partial{\cal R}_q(u^\dagger)}{\tilde{u}} \ge -\gamma_2\norm{F'(u^\dagger) \tilde{u}}
\qquad
\text{ for all } \tilde{u} \in U\;.
\end{equation}
Together with~\req{rates_sparse_q} this implies the inequality
\begin{multline}
{\cal R}_q(u)-{\cal R}_q(u^\dagger) - \inner{\partial{\cal R}_q(u^\dagger)}{u-u^\dagger}\\
\begin{aligned}\label{eq:pibot3}
&\le {\cal R}_q(u)-{\cal R}_q(u^\dagger) + \gamma_2\norm{F'(u^\dagger)(u-u^\dagger)}\\
&\le {\cal R}_q(u)-{\cal R}_q(u^\dagger) + \gamma_2 \norm{F(u)-F(u^\dagger)}\\
& \phantom{{} \le {\cal R}_q(u)-{\cal R}_q(u^\dagger)} + \gamma_2 \norm{F(u)-F(u^\dagger)-F'(u^\dagger)(u-u^\dagger)}\\
&\le (1+\gamma_2/\gamma_1)\bigl({\cal R}_q(u)-{\cal R}_q(u^\dagger)\bigr)
  + \gamma_2(1+\gamma_2/\gamma_1) \norm{F(u)-F(u^\dagger)}\;.
\end{aligned}
\end{multline}
Combination of estimates \req{pibot1}--\req{pibot3} yields
\begin{equation}\label{eq:rates_sparse_q:2}
w_{\min} \norm[b]{\pi_W^\perp u}^q \leq (1+\gamma_2/\gamma_1)\bigl({\cal R}_q(u)-{\cal R}_q(u^\dagger)\bigr)\\
+ \gamma_2(1+\gamma_2/\gamma_1) \norm{F(u)-F(u^\dagger)}\;.
\end{equation}

It remains to find an estimate for $\norm{F'(u^\dagger)(u-u^\dagger)}^q$.
Since by assumption ${\cal R}_q(u^\dagger)$, ${\cal R}_q(u) < \rho$, and $\norm{F(u)-F(u^\dagger)} < \sigma$,
it follows from~\req{rates_sparse_q} that
\begin{multline}\label{eq:rates_sparse_q:3}
\norm{F'(u^\dagger)(u-u^\dagger)}^q \\
\begin{aligned}
& \le 2^{q-1}\norm[b]{F(u)-F(u^\dagger)-F'(u^\dagger)(u-u^\dagger)}^q + 2^{q-1}\norm{F(u)-F(u^\dagger)}^q\\
& \le \frac{2^{q-1}}{\gamma_1^q} \Bigl({\cal R}_q(u)-{\cal R}_q(u^\dagger)+\gamma_2 \norm[b]{F(u)-F(u^\dagger)}\Bigr)^q
 + 2^{q-1}\norm[b]{F(u)-F(u^\dagger)}^q\\
& \le \frac{2^{q-1}(\rho+\gamma_2\sigma)^{q-1}}{\gamma_1^q}\Bigl({\cal R}_q(u)-{\cal R}_q(u^\dagger)\Bigr)\\
&\qquad\qquad\qquad
  + \Bigl(2^{q-1}\sigma^{q-1}+ \frac{2^{q-1}(\rho+\gamma_2\sigma)^{q-1}\gamma_2}{\gamma_1^q}\Bigr)\norm[b]{F(u)-F(u^\dagger)}\;.
\end{aligned}
\end{multline}
Combining the inequalities~\req{rates_sparse_q:1}, \req{rates_sparse_q:2}, and \req{rates_sparse_q:3},
we obtain the assertion.
\end{proof}

The argumentation in the proof of Theorem~\ref{th:rates_sparse_q} cannot be applied directly
to the case $q = 1$.
The main difficulty is that here the estimate~\req{subgradient}
does not follow from~\req{rates_sparse_q},
since the subgradient of ${\cal R}_1$ is not single valued.
Therefore it is necessary to postulate the existence of a subgradient element
$\xi \in \partial{\cal R}_1(u^\dagger)$ for which such an inequality holds.

\begin{theorem}[$q=1$]\label{th:rates_sparse_1}
Let $q=1$ and assume that Assumption~\ref{as:sparse_rates} holds.
In addition we assume the existence of\/ $\xi \in \partial{\cal R}_1(u^\dagger)$ and\/ $\gamma_3 > 0$ such that
\begin{equation}\label{eq:rates_sparse_1}
{\cal R}_1(u) - {\cal R}_1(u^\dagger)
\ge -\gamma_3\inner[b]{\xi}{u-u^\dagger} - \gamma_2\norm[b]{F(u)-F(u^\dagger)}
\end{equation}
for all $u \in \domain(F)$ with ${\cal R}_1(u) < \rho$ and $\norm{F(u)-F(u^\dagger)} < \sigma$.

Then it follows for a parameter choice strategy $\alpha \sim \delta^{p-1}$ that
\begin{equation*}
\norm{u_\alpha^\delta-u^\dagger} = O(\delta)\;.
\end{equation*}
\end{theorem}

\begin{proof}
We show that Assumption~\ref{as:ratesr} holds with $r = 1$. Then the result follows from
Proposition \ref{pr:con_rates}.
\smallskip

Define $J:= \set[n]{i\in\N}{\abs{\inner{\phi_i}{\xi}}\ge w_{\min}}$
and $W := \Span\set[n]{\phi_j}{j\in J}$.
Since $\xi \in U$, it follows that $J$ is a finite set.
Therefore there exists $C > 0$ such that
$C\norm{F'(u^\dagger)w} \ge \norm{w}$ for all $w \in W$.

By assumption we have that $\inner{\phi_i}{u^\dagger} = 0$ for every $i \not\in J$.
Proceeding as in the proof of Theorem~\ref{th:rates_sparse_q},
we obtain that
\[
\norm{u-u^\dagger}
\le C\norm[b]{F'(u^\dagger)(u-u^\dagger)} + \bigl(1+C\norm{F'(u^\dagger)}\bigr)\norm{\pi_W^\perp u}\;.
\]
Denote now $m := \max\set{\abs{\inner{\phi_i}{\xi}}}{i \not\in J}$,
which is well-defined, as $\bigl(\inner{\phi_i}{\xi}\bigr)_{i\in \N} \in l^2$ and therefore converges to zero.
Using the inequalities $0 \le m < w_{\min}$ and $\inner{\phi_i}{\xi} \leq m$, the assumption $\xi \in \partial{\cal R}_1(u^\dagger)$,
and \req{rates_sparse_1}, we can therefore estimate
\[
\begin{aligned}
\norm{\pi_W^\perp u}
&= \Bigl(\sum_{i\not\in J} \abs{\inner{\phi_i}{u}}^2\Bigr)^{1/2}
   \le \sum_{i\not\in J} \abs{\inner{\phi_i}{u}}\\
&\le \frac{1}{w_{\min}-m} \sum_{i\not\in J} (w_i-m)\abs{\inner{\phi_i}{u}}\\
&\le \frac{1}{w_{\min}-m} \sum_{i\not\in J} \bigl(w_i\abs{\inner{\phi_i}{u}} - \inner{\phi_i}{\xi} \inner{\phi_i}{u}\bigr)\\
&\le \frac{1}{w_{\min}-m} \sum_{i \in \N} \bigl(w_i\abs{\inner{\phi_i}{u}}
-w_i\abs{\inner{\phi_i}{u^\dagger}} - \inner{\phi_i}{\xi} \inner{\phi_i}{u-u^\dagger}\bigr)\\
&= \frac{1}{w_{\min}-m}\Bigl({\cal R}_1(u)-{\cal R}_1(u^\dagger)-\inner[b]{\xi}{u-u^\dagger}\Bigr)\\
&\le \frac{1}{w_{\min}-m}\Bigl(\bigl(1+\gamma_3^{-1}\bigr)\bigl({\cal R}_1(u)-{\cal R}_1(u^\dagger)\bigr)
   + \gamma_2/\gamma_3 \norm[b]{F(u)-F(u^\dagger)}\Bigr)\;.
\end{aligned}
\]
Here, the third to last line follows from the definition of the subgradient
and the fact that $\inner{\phi_i}{u^\dagger} = 0$ for $i \not\in J$.

For $\norm{F'(u^\dagger)(u-u^\dagger)}$ we obtain from~\req{rates_sparse_q} the estimate
\[
\begin{aligned}
\norm[b]{F'(u^\dagger)(u-u^\dagger)}
&\le \norm[b]{F(u)-F(u^\dagger)-F'(u^\dagger)(u-u^\dagger)} + \norm[b]{F(u)-F(u^\dagger)}\\
&\le \gamma_1^{-1} \bigl({\cal R}_1(u)-{\cal R}_1(u^\dagger)\bigr)
     + (1+\gamma_2/\gamma_1)\norm[b]{F(u)-F(u^\dagger)}\;.
\end{aligned}
\]
Again, the assertion follows by collecting the above inequalities.
\end{proof}

\begin{remark}
Note that in fact for the convergence rates to hold
the injectivity of $F'(u^\dagger)$ is only required
on the subspace $W$ defined in the proofs of Theorems~\ref{th:rates_sparse_q} and~\ref{th:rates_sparse_1}.
\end{remark}

\begin{remark}
Consider now the special case, where $F\colon U \to V$ is linear and bounded.
Then~\req{rates_sparse_q} with $1 \le q \le 2$ is equivalent to the source condition
\[
\partial{\cal R}_q(u^\dagger) \cap \range(F^*) \neq \emptyset\;.
\]

Indeed, in this case the operator $F$ equals its differential and therefore~\req{rates_sparse_q} reads as
\begin{equation}
\label{eq:linear_sparse}
{\cal R}_q(u) - {\cal R}_q(u^\dagger) \ge - \gamma_2 \norm[b]{F(u-u^\dagger)}\,,
\end{equation}
which is equivalent to the existence of some $\omega \in \partial{\cal R}_q(u^\dagger)$ satisfying
\[
\inner[b]{\omega}{u-u^\dagger} \ge - \gamma_2 \norm[b]{F(u-u^\dagger)}\;.
\]
This last inequality is in turn equivalent to the condition
$\partial{\cal R}_q(u^\dagger) \cap \range(F^*) \neq \emptyset$, which shows the assertion.
\smallskip

In the case $q=1$ the inequality \req{rates_sparse_1} with $\gamma_3 =1/2$ follows from \req{linear_sparse}, since
\begin{equation*}
 {\cal R}_1(u)-{\cal R}_1(u^\dagger) + \frac{1}{2} \inner[b]{\xi}{u-u^\dagger}
\ge  \frac{1}{2} \bigl({\cal R}_1(u)-{\cal R}_1(u^\dagger)\bigr)
\ge -\frac{\gamma_2}{2} \norm[b]{F(u-u^\dagger)}\;.
\end{equation*}
As a consequence, the convergence rate $O(\delta^{1/q})$ follows from the range condition
$\partial{\cal R}_q(u^\dagger) \cap \range(F^*) \neq \emptyset$ and the
finite basis injectivity property, which postulates the injectivity of the restriction
of $F$ to every subspace of $U$ spanned by a finite number of basis elements $\phi_i$.
\end{remark}

\section{Conclusion}

We have studied the application of Tikhonov regularization with $l^q$ type penalty term for
$1 \le q \le 2$ to sparse regularization.
In general, quadratic and $l^q$ regularization enjoy the same basic properties
concerning existence, stability, and convergence of the corresponding approximate solutions.
If additionally $q$ is strictly greater than one, then also the same convergence
rates can be obtained provided a source condition holds.
\smallskip

For linear operators $F$ this condition requires the subgradient  of
the penalty term to be contained in the range of the adjoint of $F$.
This assumption implies convergence rates with respect to the Bregman distance,
which for non-quadratic functionals in general cannot be compared with the norm
on the Hilbert space. In the $l^q$ case, however, such a comparison is possible
and leads to convergence rates of order $\sqrt{\delta}$ in the norm.
\smallskip

Even better results hold if the true solution $u^\dagger$ of the considered problem
is known to have a sparse representation in the chosen basis.
Then the $l^q$ regularization method yields rates of order $\delta^{1/q}$,
as long as the derivative of the operator $F$ at $u^\dagger$ is injective on the
subspace spanned by the non-zero components of $u^\dagger$.
For $q=1$ and an additional assumption concerning the subgradient of the
penalty term, this implies linear convergence of the regularized
solutions to $u^\dagger$.

\section*{Acknowledgement}

This work has been supported by the Austrian Science Fund (FWF)
within the national research networks Industrial Geometry,
project 9203-N12, and Photo\-acoustic Imaging in Biology and Medicine, project S10505-N20,
and by the Technology Transfer Office of the University of Innsbruck (transIT).

The authors want to express their thanks to Andreas Neubauer for his careful proofreading
of the article and to the referees for their valuable suggestions and comments.

\def\cprime{$'$}

\end{document}